\documentclass{amsart} \usepackage{hyperref}
\usepackage{enumerate} 
\usepackage{upref}
\usepackage{verbatim} 
\usepackage{color}
\usepackage{graphicx}
\usepackage{todonotes}
\newcommand{\sign}{\mathop{\rm sign}}
\newcommand*{\mailto}[1]{\href{mailto:#1}{\nolinkurl{#1}}}

\usepackage[latin1]{inputenc}
\usepackage{amssymb, amsmath, amsfonts, amsthm}
\usepackage[all]{xy}

\usepackage{marginnote}

\DeclareMathOperator{\id}{Id}
\DeclareMathOperator{\meas}{meas}
\DeclareMathOperator{\supp}{supp}

\newcommand{\Gr}{G}

\newcommand{\B}{\ensuremath{\mathcal{B}}}

\newcommand{\D}{\ensuremath{\mathcal{D}}}

\newcommand{\F}{\ensuremath{\mathcal{F}}}

\newcommand{\abs}[1]{\left\vert#1\right\vert}

\newcommand{\Real}{\mathbb R}

\newcommand{\norm}[1]{\left\Vert#1\right\Vert}

\newcommand{\muac}{\mu_{\text{\rm ac}}}

\newcommand{\nn}{\nonumber}

\newtheorem{theorem}{Theorem}[section]

\newtheorem{lemma}[theorem]{Lemma}
\newtheorem{definition}[theorem]{Definition}

\numberwithin{equation}{section}

\allowdisplaybreaks

\begin{document}

\title[Blow-up for the 2CH system]{Blow-up for the two-component Camassa--Holm system}

\author[K. Grunert]{Katrin Grunert}
\address{Department of Mathematical Sciences\\ Norwegian University of Science and Technology\\ NO-7491 Trondheim\\ Norway}
\email{\mailto{katring@math.ntnu.no}}
\urladdr{\url{http://www.math.ntnu.no/~katring/}}

\thanks{Research supported by the Austrian Science Fund (FWF) under Grant No.~J3147.}  

\begin{abstract}
Following conservative solutions of the two-component Camassa--Holm system $u_t-u_{txx}+3uu_x-2u_xu_{xx}-uu_{xxx}+\rho\rho_x=0$, $\rho_t+(u\rho)_x=0$ along characteristics, we determine if wave breaking occurs in the nearby future or not, for initial data $u_0\in H^1(\Real)$ and $\rho_0\in L^2(\Real)$.
\end{abstract}
\maketitle

\section{Introduction}

Over the last few years a lot of generalizations of the well-studied Camassa--Holm (CH) equation, \cite{CH:93},
\begin{equation}
 u_t-u_{txx}+3uu_{x}-2u_xu_{xx}-uu_{xxx}=0
\end{equation}
have been introduced. Among them the two-component Camassa--Holm (2CH) system,
\begin{align}
 u_t-u_{txx}+3uu_x-2u_xu_{xx}-uu_{xxx}+\rho\rho_x& =0\\ \nn 
\rho_t+(u\rho)_x&=0,
\end{align}
which we want to consider here. It has been derived as a modell for shallow water by Constantin and Ivanov \cite{CI} and has been studied intensively over the last few years, see e.g. \cite{GHR4, GHR5, GY, GL, GL2} and the references therein. 

The great interest in the CH equation and its generalizations, see e.g. \cite{ChenLiu2010, eschlechyin:07, FuQu:09, GuanKarlsenYin2010, GuanYin2010, Kuzmin, TanYin} relies on the fact that even smooth initial data might lead to classical solutions which only exist locally. In the case of the 2CH system, this is due to the fact that wave breaking can occur within finite time, that is $u_x$ becomes unbounded from below pointwise while $\norm{u}_{H^1_\Real}$  and $\norm{\rho}_{L^2_\Real}$ remain bounded, see e.g. \cite{GY,GL,GL2,GZ,Y,ZL} and the references therein. More precisely the function $u_x^2+\rho^2$ turns into a positive, finite Radon measure at breaking time, which means that energy concentrates on sets of measure zero. Thus in order to describe weak solutions, one does not only consider pairs $(u,\rho)$, such that $u\in H^1(\Real)$ and $\rho\in L^2(\Real)$, but triplets $(u,\rho,\mu)\in \D$, where $\mu$ is a positive, finite Radon measure with $\mu_{ac}=(u_x^2+\rho^2)dx$, see Definition~\ref{def:D}.

However, one is still facing the problem, what would be a natural way to prolong the solution beyond wave breaking. One possibility is to reformulate the 2CH system in Lagrangian coordinates $(y,U,h,r)\in\F$ (see Definition~\ref{th:Ldef} and \ref{def:F}). A big advantage of this approach is that the measure $\mu$ in Eulerian coordinates is mapped to the function $h$ in Lagrangian coordinates,
\newpage
\begin{subequations}
\begin{align*}
y(\xi)&=\sup\left\{y\ |\ \mu((-\infty,y))+y<\xi\right\},\\
h(\xi)&=1-y_\xi(\xi).
\end{align*}
\end{subequations}
Moreover, this new setting allows to identify where wave breaking takes place by identifying the points where $y_\xi(\xi)=0$ and in particular, all the involved functions $(y,U,h,r)$ are bounded. The associated system of differential equations \eqref{eq:chsyseq}, which describes the time evolution, provides us then with one possibility to continue the solution after wave breaking.  However, dependent on how $h$ is manipulated at breaking time, one obtains different kinds of solutions. The most prominent ones are the conservative solutions, where the energy contained remains unchanged with respect to time, see \cite{GHR4} and the dissipative ones, where a sudden drop of energy occurs at breaking time, which have been 
studied in \cite{GHR5}. Rather recently a new class of solutions has been introduced, the so-called $\alpha$-dissipative solutions, which provide a continuous interpolation between conservative and dissipative solutions, \cite{GHR6}.

The aim of this note is to show when wave breaking can occur for initial data $(u_0,\rho_0,\mu_0)\in\D$. All our considerations will be based on the description of conservative solutions in Lagrangian coordinates. Therefore we summarize in Section~\ref{sec:2} the interplay between Eulerian and Lagrangian coordinates and how the underlying system of ordinary differential equations, which describes global conservative solutions of the 2CH system, looks like. While doing so, we will focus on what wave breaking means in different formulations and what we can read off from the system of differential equations for the behavior before and after breaking time. 

Section~\ref{sec:3} then focuses on the prediction of wave breaking by following $u_x$ and $\rho$ along characteristics $y$, which solve 
\begin{equation}
y_t(t,\xi)=u(t,y(t,\xi))
\end{equation}
 for some given initial data $y_0(\xi)$. We will prove the following result.
\begin{theorem}\label{thm:main}
 Given $(u_0,\rho_0,\mu_0)\in \D$ and denote by $(u(t),\rho(t),\mu(t))\in\D$ the global conservative solution of the 2CH system at time $t$. Moreover, let $C=2(\norm{u_0}_{L^2_\Real}^2+\mu_0(\Real))$, then the following holds
\begin{enumerate}
 \item If $\rho_0(x)=0$ and $u_{0,x}(x)< -\sqrt{2C}$ for some $x\in\Real$, then wave breaking will occur within the time interval $[0,T]$, where $T$ denotes the solution of 
 \begin{equation}
 \frac{u_{0,x}(x)+\sqrt{2C}}{u_{0,x}(x)-\sqrt{2C}}= \exp(-\sqrt{2C}T).
\end{equation}
\item If $\rho_0(x)=0$ and $u_{0,x}(x)> \sqrt{2C}$ for some $x\in\Real$, then wave breaking occured within the time interval $[T,0]$, where $T$ denotes the solution of 
 \begin{equation}
 \frac{u_{0,x}(x)+\sqrt{2C}}{u_{0,x}(x)-\sqrt{2C}}= \exp(-\sqrt{2C}T).
\end{equation}
\item Assume $\rho_0(x)\not=0$ for some $x\in\Real$. Let $\xi\in\Real$ be such that $y_0(\xi)=x$, where $y_0(\xi)$ denotes the initial characteristic curve defined by \eqref{eq:Ldef1}. Denote by $y(t,\xi)$ the solution of $y_t(t,\xi)=u(t,y(t,\xi))$, then 
$u_{x}(t,y(t,\xi))$ cannot blow up.
\end{enumerate}
\end{theorem}

\section{Eulerian vs. Lagrangian coordinates}\label{sec:2}

The aim of this section is to present the interplay between Eulerian and Lagrangian coordinates and how the underlying system of ordinary differential equations in Lagrangian coordinates, which describes global conservative solutions of the 2CH system, looks like. Since the results stated in this section have already been proved earlier, we refer the interested reader for more details to \cite{GHR4}, while we here focus on pointing out properties related to wave breaking.

Let us first focus on the interplay between Eulerian and Lagrangian coordinates. It is well-known that solutions of the 2CH system might enjoy wave breaking within finite time, which means that positive amounts of energy can concentrate on sets of measure zero. 
Thus the set of possible initial data, will not consist of pairs $(u,\rho)$, where $u$ and $\rho$ are functions, but of triplets $(u,\rho,\mu)$, where $\mu$ is a positive, finite Radon measure, which describes the concentration of energy. The admissible set of Eulerian coordinates is then given by
\begin{definition} [Eulerian coordinates]\label{def:D} The set $\D$ is
  composed of all triplets $(u,\rho,\mu)$ such that
  $u\in H^1(\Real)$, $\rho\in L^2(\Real)$ and $\mu$ is a positive, finite
  Radon measure whose absolutely continuous part,
  $\muac$, satisfies
\begin{equation}
\label{eq:abspart}
\muac=(u_x^2+\rho^2)\,dx.
\end{equation}
\end{definition}

The weak formulation of the 2CH system
\begin{subequations}\label{equiv:sys}
\begin{align}
 u_t+uu_{x}+P_x&=0,\\
\rho_t+(u\rho)_x& =0,
\end{align}
\end{subequations}
where 
\begin{equation}
P(t,x)=\frac12 \int_\Real e^{-\vert x-z\vert}u^2(t,z)dz+\frac14 \int_\Real e^{-\vert x-z\vert}d\mu(t,z),
\end{equation}
hints that it could be possible to describe solutions of the 2CH system with the help of the method of characteristics. Indeed, one defines the time evolution of the characteristic $y(t,\xi)$ to be 
\begin{equation}\label{eq:char}
 y_t(t,\xi)=u(t,y(t,\xi)),
\end{equation}
for some given initial data $y_0(\xi)=y(0,\xi)$,
and introduces the functions 
\begin{equation}\label{eq:map}
 U(t,\xi)=u(t,y(t,\xi)) \quad \text{ and }\quad r(t,\xi)=\rho(t,y(t,\xi))y_\xi(t,\xi), 
\end{equation}
whose time evolution is given through \eqref{equiv:sys}. However, this ansatz has two drawbacks, on the one hand we do not know yet how $\mu(t,x)$ evolves with respect to time, and on the other hand how to combine $\mu(t,x)$ and $y(t,\xi)$. The last issue will be resolved next by introducing a suitable mapping from Eulerian coordinates $\D$ to Lagrangian coordinates $\F$.

\begin{definition} 
\label{th:Ldef}
For any $(u,\rho,\mu)$ in $\D$, let
\begin{subequations}
\label{eq:Ldef}
\begin{align}
\label{eq:Ldef1}
y(\xi)&=\sup\left\{y\ |\ \mu((-\infty,y))+y<\xi\right\},\\
\label{eq:Ldef2}
h(\xi)&=1-y_\xi(\xi),\\
\label{eq:Ldef3}
U(\xi)&=u\circ{y(\xi)},\\
\label{eq:Ldef4}
r(\xi)&=\rho\circ{y(\xi)}y_\xi(\xi).
\end{align}
\end{subequations}
Then $(y,U,h,r)\in\F$. We denote by $L\colon \D\rightarrow \F$ the mapping which to any element $(u,\rho,\mu)\in\D$ associates $X=(y,U,h,r)\in \F$ given by \eqref{eq:Ldef}.  
\end{definition}

Since $\mu$ is a positive, finite Radon measure, we have that the function $F(x)=\mu((-\infty,x))$ is increasing and lower semi-continuous. We therefore obtain, as an immediate consequence, that $y$ is an increasing function, such that we can find to any $x\in\Real$ at least one $\xi\in\Real$ such that $y(\xi)=x$. In particular, we have that if a positive amount of energy concentrates at some point $x$, i.e. $\mu(\{x\})=c>0$, then \eqref{eq:Ldef1} implies that the point $x$ will be mapped to an interval $[\xi_1,\xi_2]$, such that $\xi_2-\xi_1=c$ and $y(\xi)=x$ for all $\xi\in[\xi_1,\xi_2]$. Moreover, $y$ is Lipschitz continuous with Lipschitz constant at most one, and hence differentiable almost everywhere by Rademacher's theorem. Thus $y_\xi(\xi)=0$ for almost every $\xi\in[\xi_1,\xi_2]$, which is equivalent to $h(\xi)=1$ for almost every $\xi\in[\xi_1,\xi_2]$. Thus $\mu(\{x\})=\int_{\xi_1}^{\xi_2}h(\xi)d\xi$. In general, one has that 
\begin{equation*}
\{\xi\in\Real \mid y_\xi(\xi)=0\}=\{\xi\in\Real\mid y(\xi)\in\supp(\mu_s)\},
\end{equation*}
by Besicovitch's derivation theorem, see e.g. \cite[Theorem 2.22]{Ambrosio}, and $\mu(\Real)=\norm{h}_{L^1_\Real}$. 

The above considerations show that wave breaking in Lagrangian coordinates means that $y_\xi(\xi)$ equals $0$, and the concentrated energy can be read of from the function $h(\xi)$ integrated over $\{\xi\in\Real\mid y_\xi(\xi)=0\}$. In particular, one can view, roughly speaking, $h$, which is a function, as the relabeled version of the measure $\mu$. Therefore, \eqref{eq:abspart} implies that $hy_\xi=U_\xi^2+r^2$ for almost every $\xi\in\Real$. Since $y_\xi$ and $h$ are bounded, we have that $y_\xi(\xi)=0$ implies $U_\xi(\xi)=r(\xi)=0$.

Before introducing the set of Lagrangian coordinates, which provides us with a list of all properties inherited by $(y,U,h,r)$ necessary for global conservative solutions to exist, we will characterize relabeling functions. These functions will not only play a key role when defining $\F$, but also when identifying equivalence classes in Lagrangian coordinates, which turn up since we introduced 4 Lagrangian coordinates in contrast to 3 Eulerian coordinates, and hence a degree of freedom, which could  cause trouble, when considering questions concerning the uniqueness of solutions. 

\begin{definition}[Relabeling functions]\label{def:rf}
 We denote by $G$ the subgroup of the group of homeomorphisms from $\Real$ to $\Real$ such that 
\begin{subequations}
\label{eq:Gcond}
 \begin{align}
  \label{eq:Gcond1}
  f-\id \text{ and } f^{-1}-\id &\text{ both belong to } W^{1,\infty}(\Real), \\
  \label{eq:Gcond2}
  f_\xi-1 &\text{ belongs to } L^2(\Real),
 \end{align}
\end{subequations}
where $\id$ denotes the identity function. Given $\kappa>0$, we denote by $G_\kappa$ the subset $G_\kappa$ of $G$ defined by 
\begin{equation}
 G_\kappa=\{ f\in G\mid  \norm{f-\id}_{W^{1,\infty}}+\norm{f^{-1}-\id}_{W^{1,\infty}}\leq\kappa\}. 
\end{equation}
\end{definition}

At first sight it seems difficult to check if a function belongs to $G$ or not, the following lemma simplifies this task considerably.
\begin{lemma}[{\cite[Lemma 3.2]{HR}}] 
\label{lem:charH}
Let $\kappa\geq0$. If $f$ belongs to $G_\kappa$,
then $1/(1+\kappa)\leq f_\xi\leq 1+\kappa$ almost
everywhere. Conversely, if $f$ is absolutely
continuous, $f-\id\in W^{1,\infty}(\Real)$, $f$ satisfies
\eqref{eq:Gcond2} and there exists $d\geq 1$ such
that $1/d\leq f_\xi\leq d$ almost everywhere, then
$f\in G_\kappa$ for some $\kappa$ depending only
on $d$ and $\norm{f-\id}_{W^{1,\infty}}$.
\end{lemma}

We are now ready to introduce the set of Lagrangian coordinates, which provides us with a list of all properties of $(y,U,h,r)$, which are necessary for global solutions to exist.
\begin{definition}[Lagrangian coordinates]\label{def:F} The set $\F$ is composed of all $X=(\zeta,U,h,r)$, such that 
 \begin{subequations}
\label{eq:lagcoord}
\begin{align}
\label{eq:lagcoord1}
& (\zeta, U,h,r,\zeta_\xi,U_\xi)\in L^\infty(\Real)\times [L^\infty(\Real)\cap L^2(\Real)]^5,\\
\label{eq:lagcoord2}
&y_\xi\geq0, h\geq0, y_\xi+h>0
\text{  almost everywhere},\\
\label{eq:lagcoord3}
&y_\xi h=U_\xi^2+r^2\text{ almost everywhere},\\
\label{eq:lagcoord4}
&y+H\in G,
\end{align}
\end{subequations}
where we denote $y(\xi)=\zeta(\xi)+\xi$ and $H(\xi)=\int_{-\infty}^\xi h(\eta)d\eta$.
\end{definition}

Note that \eqref{eq:lagcoord3}, together with \eqref{eq:char} and \eqref{eq:map}, provides us with a possibility to define $h_t(t,\xi)$ for all $\xi$ such that $y_\xi(t,\xi)\not=0$. Since we want to obtain solutions, which are continuous with respect to time, the time evolution of the 2CH system in Lagrangian coordinates is given through
\begin{subequations}
  \label{eq:chsyseq}
  \begin{align}
    \label{eq:chsyseq1}
    \zeta_t&=U,\\
    \label{eq:chsyseq:2}
    \zeta_{\xi,t}& =U_\xi,\\
    \label{eq:chsyseq3}
    U_t&=-Q,\\
    \label{eq:chsyseq4}
    U_{\xi,t}& =\frac{1}{2}h+(U^2-P)y_\xi, \\
    \label{eq:chsyseq5}
    h_t&=2(U^2-P)U_\xi,\\
    \label{eq:chsyseq6}
    r_t&=0,
  \end{align}
\end{subequations}
where 
\begin{equation}
 \label{eq:P}
P(t,\xi)=\frac14\int_\Real
  e^{-\abs{y(t,\xi)-y(t,\eta)}}(2U^2y_\xi+h)(t,\eta)\,d\eta,
\end{equation}
and 
\begin{equation}
 \label{eq:Q}
Q(t,\xi)=-\frac14\int_\Real \sign(\xi-\eta)e^{-\abs{y(t,\xi)-y(t,\eta)}}(2
  U^2y_\xi+h)(t,\eta)\,d\eta.
\end{equation}
One can show that both $P$ and $Q$ belong to $H^1(\Real)$.

This system, which yields for any initial data $X_0\in\F$ a unique global solution in $\F$, describes the conservative solutions. By conservative we mean that the energy, which is given in the case of the 2CH system through 
\begin{equation}\label{eq:energy}
 E(t)=\int_\Real (U^2y_\xi+h)(t,\xi)d\xi=\norm{u(t,.)}_{L^2_\Real}^2+\mu_t(\Real),
\end{equation}
is preserved for all times, i.e. $E(t)=E(0)$ for all $t\in\Real$. Note that $E(t)$ is well-defined for all $(y,U,h,r)\in \F$. Moerover, it is not only possible to consider \eqref{eq:chsyseq} forward in time, but also backward. 

We already saw that wave breaking at time $t_b$ for some $\xi\in\Real$ means that $y_\xi(t_b,\xi)=U_\xi(t_b,\xi)=r(t_b,\xi)=0$. Then we get from \eqref{eq:chsyseq} that $h_t(t_b,\xi)=y_{\xi,t}(t_b,\xi)=r_t(t_b,\xi)=0$, while $U_{\xi,t}(t_b,\xi)=\frac12 h(t_b,\xi)$. Applying Lemma~\ref{lem:charH} to $y+H\in G$, we get that $h(t_b,\xi)=h(t_b,\xi)+y_\xi(t_b,\xi)>0$,  which means that $U_\xi(t,\xi)$ changes from negative to positive at breaking time $t_b$. We will see some consequences of this observation in the next section.

Although relabeling will not play an important role in our further investigations, we will state the most important results here, for the sake of completeness. For any $X=(y,U,h,r)\in\F$ and any function $f\in\Gr$ we denote $(y\circ f, U\circ f, h\circ ff_\xi, r\circ ff_\xi)$ by $X\circ f$. We say that $X$ and $X'\in\F$ are equivalent, if there exists a relabeling function $f\in\Gr$ such that $X'=X\circ f$. 
Let 
\begin{equation}
 \F_0=\{ X\in\F\mid y+H=\id\},
\end{equation}
then $\F_0$ contains exactly one element of each equivalence class. In particular, we have that $L$ is a mapping from $\D$ to $\F_0$ and therefore different elements in $\D$ are mapped to different elements in $\F_0$, which is why the mapping $L$ is well-defined. Moerover, let us denote by $S_t(X)$, the solutions of \eqref{eq:chsyseq} at time $t$ with initial data $X$, then 
\begin{equation}
 S_t(X\circ f)=S_t(X)\circ f.
\end{equation}
Since $X\in \F_0$ does not imply $S_t(X)\in \F_0$ for all $t\in\Real$, it would be a huge advantage if we could find a mapping from Lagrangian coordinates $\F$ to Eulerian coordinates $\D$, such that any two elements in $\F$ belonging to the same equivalence class are mapped to the same element in Eulerian coordinates.

\begin{theorem}
\label{th:umudef} 
Given any element $X=(y,U,h,r)\in\F$. Then, the
measure $y_{\#}(r(\xi)\,d\xi)$ is absolutely
continuous, and we define $(u,\rho,\mu)$ as follows
\begin{subequations}
\label{eq:umudef}
\begin{align}
\label{eq:umudef1}
u(x)&=U(\xi) \text{ for any }\xi \text{ such that }x=y(\xi),\\
\label{eq:umudef2}
\mu&=y_\#(h(\xi)\,d\xi),\\
\label{eq:umudef3}
\rho(x)\,dx&=y_\#(r(\xi)\,d\xi),
\end{align}
\end{subequations}
We have that $(u,\rho,\mu)$ belongs to $\D$. We
denote by $M\colon \F\rightarrow\D$ the mapping
which to any $X$ in $\F$ associates the element $(u, \rho,\mu)\in \D$ as given
by \eqref{eq:umudef}. In particular, the mapping $M$ is
invariant under relabeling.
\end{theorem}

Finally we recall that the push-forward of a measure $\nu$ by a measureable function $f$ is the measure $f_\#\nu(B)=\nu(f^{-1}(B))$. Thus Theorem~\ref{th:umudef} implies that if there exists an interval of positive length $[\xi_1,\xi_2]$ such that $y(\xi)=x$ for all $\xi\in[\xi_1,\xi_2]$, then $\mu(\{x\})>0$ and $\rho(x)=0$. 

\section{Wave breaking or not?}\label{sec:3}

After recalling the interplay between Eulerian and Lagrangian coordinates and presenting the system of differential equations describing the time evolution in Lagrangian coordinates, we want to investigate when wave breaking might occur, that means when $u_x$ can become unbounded from below pointwise. 

\begin{proof}[Proof of Theorem 1.1]
Given some initial data $(u_0,\rho_0,\mu_0)\in \D$, we can find to almost every $x\in\Real$, a $\xi\in\Real$ such that $x=y(0,\xi)$ and $y_\xi(0,\xi)>0$. Let us denote 
\begin{equation}
 \B=\{ \xi\in\Real\mid y(0,\xi), U(0,\xi)\text{ are differentiable and }y_\xi(0,\xi)>0\}.
\end{equation}
Then the set $y(\B)$ has full measure. Indeed, we have, after a change of variables 
\begin{equation}
 \meas (y(\B)^c)=\int_{y(\B)^c} dx=\int_{\B^c}y_\xi d\xi=0.
\end{equation}
Then for all $\xi\in\B$, the variables 
\begin{equation}
 \alpha(0,\xi)=\frac{U_\xi}{y_\xi}(0,\xi)\quad \text{ and }\quad \beta(0,\xi)=\frac{r}{y_\xi}(0,\xi),
\end{equation}
are well-defined.
Since $U_\xi=u_x\circ yy_\xi$ and $r=\rho\circ yy_\xi$, we have 
\begin{equation}
 u_x(0,x)=u_x(0,y(0,\xi))=\frac{U_\xi(0,\xi)}{y_\xi(0,\xi)}=\alpha(0,\xi),
\end{equation}
and 
\begin{equation}
 \rho(0,x)=\rho(0,y(0,\xi))=\frac{r(0,\xi)}{y_\xi(0,\xi)}=\beta(0,\xi).
\end{equation}
From \eqref{eq:chsyseq}, we get 
\begin{equation*}
 \alpha_t=\frac{U_{\xi,t}y_\xi-U_\xi y_{\xi,t}}{y_\xi^2} 
 =\frac{\frac12 hy_\xi+(U^2-P)y_\xi^2-U_\xi^2}{y_\xi^2}
 =\frac12\beta^2-\frac12 \alpha^2+(U^2-P)
\end{equation*}
and 
\begin{equation*}
 \beta_t =-\frac{ry_{\xi,t}}{y_\xi^2}=- \frac{rU_\xi}{y_\xi^2}=-\alpha\beta.
\end{equation*}
Thus the system under consideration will be 
\begin{subequations}\label{eq:sysnew}
\begin{align}\label{eq:11}
 \alpha_t&=\frac12 \beta^2-\frac12 \alpha^2+(U^2-P),\\ \label{eq:22}
\beta_t&=-\alpha\beta,
\end{align}
\end{subequations}
and we want to find out when $\alpha$ can become unbounded from below, or in other words when wave breaking occurs. 

Here it is important to note, that for $u_0\in H^1(\Real)$ and $\rho_0\in L^2(\Real)$, the energy $\norm{u(t)}_{L^2_\Real}+\mu_t(\Real)$, is preserved (cf. \eqref{eq:energy}), since we are considering conservative solutions. Thus $P$ and $P-U^2$ can be bounded uniformly (with respect to time) by a constant only dependent on $\norm{u_0}^2_{L^2_\Real}$ and $\mu_0(\Real)$. In particular, we have 
\begin{equation}
 \norm{P(t,.)}_{L^\infty_\Real}\leq \norm{u_0}_{L^2_\Real}^2+\mu_0(\Real),
\end{equation}
and 
\begin{equation}
 \norm{U^2(t,.)-P(t,.)}_{L^\infty_\Real}\leq 2(\norm{u_0}_{L^2_\Real}^2+\mu_0(\Real)).
\end{equation}

\begin{figure}\label{Interaction}
 \centering
      \includegraphics[width=.4\textwidth]{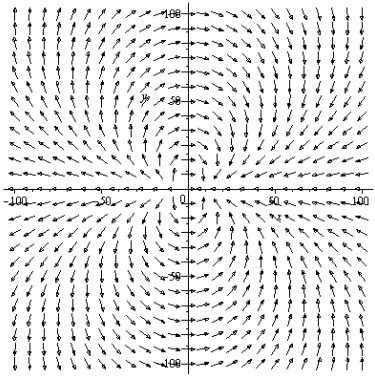}
      \caption{Plot of the vectorfield for the functions $(\alpha(t),\beta(t))$ (with $(U^2-P)$ replaced by the constant value $5$).}
\end{figure}

We will distinguish two cases for fixed $\xi\in\B$:

(\textit{i}): If $\beta(0,\xi)=0$, \eqref{eq:22} implies that $\beta(t,\xi)=0$ for all $t\geq 0$ and \eqref{eq:11} reduces to
\begin{equation}\label{eq:1}
\alpha_t= -\frac12 \alpha^2+(U^2-P).
\end{equation}
Thus we see that $\alpha_t(t,\xi)\leq 0$ for $\alpha^2(t,\xi)\geq 4(\norm{u_0}_{L^2_\Real}^2+\mu_0(\Real))$.
 
Assume that $\alpha(0,\xi)\leq -2\sqrt{\norm{u_0}_{L^2_\Real}^2+\mu_0(\Real)}$, then we can find to any constant $C<\alpha(0,\xi)<-2\sqrt{\norm{u_0}_{L^2_\Real}^2+\mu_0(\Real)}$, a time $t_C$, such that $\alpha(t_C,\xi)=C$, and $\alpha(t,\xi)$ is strictly decreasing on the interval $[0,t_C]$. 

Let $C_n$ be a decreasing sequence of constants such that  $C_n<\alpha(0,\xi)$ for all $n\in\mathbb{N}$ and $C_n\to -\infty$. Then we can associate to each $C_n$, as before, a time $t_{n}$, such that $\alpha(t_n,\xi)=C_n$ and $\alpha(t,\xi)$ is strictly decreasing on the interval $[0,t_n]$. Thus 
\begin{equation}
 \alpha(t_n,\xi)=C_n\to -\infty, \quad \text{ as } \quad n\to\infty. 
\end{equation}

Unfortunately, this argument has one drawback so far, we do not know if the sequence $t_n$ is convergent. However what we do know, is that the sequence $t_n$ is increasing. This means if we can show that the sequence $t_n$ is uniformly bounded then there exists a positive time $t_b$ such that 
\begin{equation}\label{conv:tn}
 t_n\uparrow t_b, \quad \text{ as }\quad n\to\infty,
\end{equation}
and $t_b$ is the time where wave breaking occurs, since $\alpha(t,\xi)$ becomes unbounded from below.

Therefore consider the differential equation 
\begin{equation}\label{eq:gamma}
 \gamma_t(t)=-\frac12 \gamma^2(t)+C,
\end{equation}
where $C=2(\norm{u_0}_{L^2_\Real}^2+\mu_0(\Real))$, with initial data $\gamma(0)=\alpha(0,\xi)$.  Then $\gamma^2(0)-2C> 0$ and $\gamma(t)$ is strictly decreasing. Moreover, $\gamma(t)$ is negative and an upper bound for $\alpha(t,\xi)$, i.e. $\alpha^2(t,\xi)\geq \gamma^2(t)>2C$, as long as $\alpha(t,\xi)$ is well-defined. A qualitative analysis of \eqref{eq:gamma} shows that $\gamma(t)$, will blow up. We claim that this happens within finite time. If so the blow up time of $\gamma(t)$ provides us with an upper bound for the sequence $t_n$. 
 
The solution of \eqref{eq:gamma} is given by
\begin{equation}
 \gamma(t)=\frac{\sqrt{2C}\gamma(0)+2C+(\sqrt{2C}\gamma(0)-2C)\exp(-\sqrt{2C}t)}{\gamma(0)+\sqrt{2C}-(\gamma(0)-\sqrt{2C})\exp(-\sqrt{2C}t)}.
\end{equation}
 A close look reveals that both the numerator, which remains negative for all times, and the denominator, which is positive at initial time, are absolutely bounded. Thus $\gamma^2(t)$ can only blow up within finite time if the denominator equals zero, i.e. 
\begin{equation}\label{est:t}
 \frac{\gamma(0)+\sqrt{2C}}{\gamma(0)-\sqrt{2C}}= \exp(-\sqrt{2C}t).
\end{equation}
Solving this very last equation for $t$, provides us with an upper bound on the sequence $t_n$ and proves \eqref{conv:tn}.

Although our argument is rigorous for a particular sequence $C_n$, it is left to show that $t_b$ is independent of the sequence $C_n$ we choose. Assume the opposite, that is there exist two sequences $C_n$ and $\tilde C_n$ as above, with 
\begin{equation}
\lim_{n\to\infty}t_n=t_b<\tilde t_b=\lim_{n\to\infty}\tilde t_n.
\end{equation}
Then we can find $N\in \mathbb{N}$, such that $t_n < \tilde t_N$ for all $n\in \mathbb{N}$. Since we have, by construction, that $\alpha(\tilde t_N,\xi)=\tilde C_N$, which means in particular finite, and $\alpha(t,\xi)$ is strictly decreasing on $[0,\tilde t_N]$, it follows that $C_n\geq \tilde C_N$ for all $n\in\mathbb{N}$, which contradicts $C_n\to -\infty$ as $n\to\infty$.

As already pointed out, a rigorous analysis of \eqref{eq:chsyseq} reveals that $U_\xi(t,\xi)$ changes sign from negative to positive, at breaking time.
Due to the conditions one has to impose, to obtain global conservative solutions, we have $y_\xi(t,\xi)\geq 0$ for almost every $\xi\in\Real$, which implies that $\alpha(t,\xi)$ becomes positive after wave breaking occured. 

Since all our considerations are based on the description of global conservative solutions, it is not only possible to follow $u_x$ and $\rho$ along characteristics forward in time, but also backwards. Thus a natural question that arises in this context, is, if we can find out how a solution behaves after wave breaking. 

Following the same lines as before one can show that 
 \begin{equation}
  \alpha(t,\xi)\to \infty, \quad \text{ as }\quad t\downarrow t_b.
 \end{equation}
To be a lit more explicit, one looks at the differential equation \eqref{eq:1} backwards in time for some initial data $\alpha(0,\xi)$  with $\alpha(0,\xi)\geq 2\sqrt{\norm{u_0}_{L^2_\Real}^2+\mu_0(\Real)}$ and compares it with the solution of \eqref{eq:gamma} with $\gamma(0)=\alpha(0,\xi)$ backwards in time. One then obtains that $2\sqrt{\norm{u_0}_{L^2_\Real}^2+\mu_0(\Real)}\leq\gamma(t)\leq \alpha(t,\xi)$ as long as $\alpha(t,\xi)$ is well-defined. Finally, $\gamma(t)$ will provide us with an upper bound on how much time has past since wave breaking occured for the points under consideration.

Last but not least there is one important remark. We assumed throughout all considerations that wave breaking means that $u_x$ becomes unbounded from below along characteristics, but why can $u_x$ not tend to $\infty$ along characteristics forward in time? \eqref{eq:1} implies that $\alpha(t,\xi)\leq \max \big(\alpha(0,\xi), 2\sqrt{\norm{u_0}_{L^2_\Real}^2+\mu_0(\Real)}\big)$ as long as $\alpha(t,\xi)$ is well-defined. 

(\textit{ii}): If $\beta(0,\xi)>0$ (the case $\beta(0,\xi)<0$ is similar), we get from \eqref{eq:22} that 
\begin{equation}\label{sol:beta}
\beta(t,\xi)=\beta(0,\xi)\exp(-\int_0^t \alpha(s,\xi)ds),
\end{equation}
and hence $\beta(t,\xi)\geq 0 $ for all $t\geq 0$. In particular, we see that $\beta(t,\xi)$ increases if $\alpha(t,\xi)$ is negative and $\beta(t,\xi)$ decreases if $\alpha(t,\xi)$ is positive. Note additionally that $\beta(t,\xi)$ remains finite on some interval $[0,T]$, if $\alpha(t,\xi)$ is bounded from below on $[0,T]$.

We now want to show that no wave breaking can occur in this case. Let us assume the opposite. Namely, $\alpha(t,\xi)$ becomes unbounded from below within finite time for some fixed $\xi\in\Real$. Denote by $t_b$ the first time when wave breaking occurs, then 
\begin{equation}
 \alpha(t,\xi)\to -\infty \quad \text{ as }\quad t\to t_b. 
\end{equation}
Moreover, \eqref{sol:beta} implies that either $\beta(t_b,\xi)$ is  positive and finite or 
\begin{equation}
 \beta(t,\xi)\to \infty \text{ as }t\to t_b.
\end{equation}

Since the right hans side of \eqref{eq:sysnew} depends on $\alpha^2(t,\xi)$, $\beta^2(t,\xi)$ and the uniformly bounded function $U^2(t,\xi)-P(t,\xi)$,
we are going to investigate the interplay between $\alpha^2(t,\xi)$ and $\beta^2(t,\xi)$. Therefore consider the function
\begin{equation}
\frac{\alpha^2(t,\xi)}{\beta^2(t,\xi)},
\end{equation}
which is bounded from below by $0$. According to \eqref{eq:sysnew}, we have 
\begin{equation}\label{eq:33}
 \Big(\frac{\alpha^2}{\beta^2}\Big)_t=\frac{\alpha}{\beta^2}(\alpha^2+\beta^2+2(U^2-P)).
\end{equation}
Since $\alpha(t,\xi)\to -\infty$ as $t\to t_b$, there exists a time $0\leq \tilde t\leq t_b$, such that $\alpha_t(\tilde t,\xi)<0$, that is
\begin{equation}\label{eq:beg}
 \alpha^2(\tilde t,\xi)\geq \beta^2(\tilde t,\xi)+2(U^2(\tilde t,\xi)-P(\tilde t,\xi)),
\end{equation}
and
\begin{equation}\label{eq:end}
 \alpha^2(t,\xi)\geq 8 (\norm{u_0}_{L^2_\Real}^2+\mu_0(\Real))\geq 4\norm{U^2(t,.)-P(t,.)}_{L^\infty_\Real},
\end{equation}
for all $t\in [\tilde t, t_b)$. Note that the last inequality implies that $\alpha(t,\xi)<0$ for all $t\in [\tilde t, t_b)$.
Furthermore, $\frac{\alpha^2(\tilde t, \xi)}{\beta^2(\tilde t,\xi)}$ is finite and strictly positive due to \eqref{eq:beg} and \eqref{eq:end}. Combining now \eqref{eq:33} and \eqref{eq:end}, yields that the function $\frac{\alpha^2(t,\xi)}{\beta^2(t,\xi)}$ is decreasing on the interval $[\tilde t, t_b)$, since $\alpha(t,\xi)$ is negative for all $t\in[\tilde t, t_b)$. 
This implies, since 
\begin{equation}
0\leq \frac{\alpha^2(t,\xi)}{\beta^2(t,\xi)}\leq \frac{\alpha^2(\tilde t, \xi)}{\beta^2(\tilde t, \xi)} \quad \text{ for all }\quad t\in[\tilde t, t_b),
\end{equation}
that $\lim_{t\uparrow t_b} \frac{\alpha^2(t,\xi)}{\beta^2(t,\xi)}$ exists. Additionally, we have by \eqref{eq:33} and \eqref{eq:end}, that 
\begin{align}
 0\leq \int_{t_b}^{\tilde t} \alpha(s,\xi)ds
&\leq\int_{t_b}^{\tilde t} \Big(\alpha(s,\xi)+\alpha(t,\xi)\Big(\frac{\alpha^2(s,\xi)}{\beta^2(s,\xi)}+\frac{2(U^2(s,\xi)-P(s,\xi))}{\beta^2(s,\xi)}\Big)\Big)ds\\ \nn
& = \frac{\alpha^2(\tilde t,\xi)}{\beta^2(\tilde t,\xi)}-\frac{\alpha^2(t_b,\xi)}{\beta^2(t_b,\xi)}<\infty.
\end{align}
Thus, according to \eqref{sol:beta}, $\lim_{t\uparrow t_b}\beta(t,\xi)$ exists and 
\begin{equation}
 \lim_{t\uparrow t_b} \beta(t,\xi)=\beta(\tilde t,\xi)\exp\Big(-\int_{\tilde t}^{t_b} \alpha(s,\xi)ds\Big).
\end{equation}
Moerover, we get 
\begin{equation}
 \lim_{t\uparrow t_b} \alpha^2(t,\xi)=\lim_{t\uparrow t_b} \beta^2(t,\xi)\lim_{t\uparrow t_b}\frac{\alpha^2(t,\xi)}{\beta^2(t,\xi)}=\beta^2(t_b,\xi)\frac{\alpha^2(t_b,\xi)}{\beta^2(t_b,\xi)}<\infty,
\end{equation}
and in particular $\alpha(t_b,\xi)$ is finite and no blow up occurs, which contradicts our assumption.

Note, that the above considerations imply that if $\alpha_t(t,\xi)<0$, then $\beta^2(t,\xi)$ is growing faster than $\alpha^2(t,\xi)$. This is the main reason why wave breaking cannot occur.

\end{proof}

Last but not least we want to point out that if $\beta(0,\xi)\not=0$, then $u_x$ cannot tend to $\infty$ along characteristics. The proof follows the same lines as the one of Theorem~\ref{thm:main} (\textit{iii}), with slight modifications.

\end{document}